\documentclass[12pt,psamsfonts]{amsart}
 \usepackage{amssymb}
\usepackage[colorlinks=true, pdfstartview=FitV, linkcolor=blue, citecolor=blue, urlcolor=blue]{hyperref}
\usepackage{amscd} 
\usepackage{longtable}
\calclayout

\numberwithin{equation}{subsection}

\newtheorem{theorem}[subsection]{Theorem}
\newtheorem{lemma}[subsection]{Lemma}
\newtheorem{proposition}[subsection]{Proposition}
\newtheorem{corollary}[subsection]{Corollary}

\theoremstyle{definition}

\theoremstyle{remark} \newtheorem{remark}[subsection]{Remark}

\newtheorem{example}[subsection]{Example}

\def\inv{^{-1}} 
\def\Ker{\operatorname{Ker}}
\def\today{\ifcase\month\or January\or
  February\or March\or April\or May\or June\or July\or August\or
  September\or October\or November\or December\fi \space\number\day,
  \number\year} 
  \def\C{\mathbb C} 
  \def\O{\mathcal O}
\def\twist#1#2#3{#1\ast^{#2}#3}
  \def\lie#1{\mathfrak{ #1}}
\def\lieh{\lie h} 
\def\liek{\lie k} 
\def\lieh{\lie h} 
\def\lieg{\lie g} 
\def\lies{\lie s}
\def\R{{\mathbb R}}

\def\pt{\partial} 
\def\SU{\operatorname{SU}}

 \def\GL{\operatorname{GL}}
 
\def\Aut{\operatorname{Aut}}
\def\X{ {\mathcal X}}
\def\H{{\mathcal H}}
\def\xck{\X^\infty_c(X)^K}
\def\xk{\X^\infty(X)^K}
\def\supp{\operatorname{supp}}
\def\cc{\mathcal C}
\def\End{\operatorname{End}}
\def\liegl{{\lie{gl}}}
\def\HH{\mathbb H}
\def\Der{\operatorname{Der}}
\def\Im{\operatorname{Im}}
\def\B{\mathcal{B}}
\def\phi{\varphi}
\def\Id{\operatorname{Id}}
\title{Homology of equivariant vector fields}   \author{Gerald W.
  Schwarz} \thanks{Partially supported by NSA grant
  H98230--04--01--0070.}
\address{Department of Mathematics\\
  Brandeis University\\
  PO Box 549110\\
  Waltham, MA 02454-9110} \email{schwarz@brandeis.edu}

\begin{document}
 \begin{abstract}
 Let $K$ be a compact Lie group. We compute the abelianization of the Lie algebra of equivariant vector fields on a smooth $K$-manifold $X$. We also compute the abelianization of the Lie algebra of strata preserving smooth vector fields on the quotient $X/K$.
 \end{abstract}
\maketitle

 \section{Introduction}
\subsection{} K.~Abe and K.~Fukui \cite{AbeFukui2} have considered the first homology group (abelianization) of the group  of equivariant smooth diffeomorphisms of a smooth $K$-manifold $X$, where $K$ is finite. They also computed the abelianization for the diffeomorphisms of the quotient orbifold $X/K$. Our results below are the analogues of their results for  vector fields in the case that $K$ is a compact Lie group. The  vector fields are, in  a sense, the Lie algebras of the relevant diffeomorphism groups, so, hopefully, our results  indicate that one should be able to generalize   the Abe-Fukui results. There are already generalizations in some cases  \cite{AbeFukui1}.

\subsection{} Let $X$ be a smooth $K$-manifold where $K$ is compact. Let $\X^\infty(X)$ denote the Lie algebra of smooth vector fields on $X$ and let $\X^\infty_c(X)$ denote the subalgebra of vector fields with compact support. If $X$ is algebraic, then $\X(X)$ will denote the polynomial vector fields on $X$. By $\X^\infty(X)^K$, etc.\ we mean the $K$-invariant elements in $\X^\infty(X)$, etc. We will state most of our results  for $\X^\infty_c(X)^K$; the corresponding results for $\X^\infty(X)^K$ follow easily from our techniques.
 
 If $\lieg$ is a Lie algebra, we denote by $\H(\lieg)$ the abelianization $\lieg/[\lieg,\lieg]$. We denote the Lie algebras of   compact Lie groups $K$, $H$, etc.\ by the corresponding gothic letters $\liek$, $\lieh$, etc.
 
\subsection{} Let $x\in X$. Then we have the isotropy group $K_x$ and its slice representation on $W_x:=T_xX/T_x(Kx)$ where $Kx$ denotes the $K$-orbit through $x$. We say that the orbit $Kx$ is {\it  isolated\/} if $W_x^{K_x}=(0)$. It follows from the differentiable slice theorem that $Kx$ is isolated if and only if all  isotropy groups $K_y$ of   points $y$ near $x$, $Ky\neq Kx$, are conjugate to a proper subgroup of $K_x$.  There is then a discrete subset $\{x_i\}_{i\in I}$ of $X$ (possibly empty) where we choose one point from each isolated orbit. Let $H_i$ denote $K_{x_i}$ and set $W_i:=W_{x_i}$, $i\in I$. 
 
 \begin{theorem}\label{maintheorem}
  Let $X$  and the $x_i$, $H_i$ and $W_i$ be as above Then
  $$
  \H(\X^\infty_c(X)^K)\simeq\bigoplus_i\H(\liek^{H_i}/\lieh_i^{H_i})\bigoplus_i\H(\End(W_i)^{H_i}).
  $$
  \end{theorem}
   \begin{theorem}\label{localtheorem} Let $H$ be a compact Lie group and $V$ an $H$-module where $V^H=(0)$.
   Write $V=\oplus_{j=1}^m n_j V_j$ where the $V_j$ are irreducible and pairwise non-isomorphic and $n_jV_j$ denotes the direct sum of $n_j$ copies of $V_j$. Let  $l$ denote the number of $V_j$ such that $\End(V_j)^{H}\simeq\C$ and let $Z(\End(V)^H)$ denote the center of $\End(V)^H$. Then
  $$
  \H(\End(V)^H)\simeq Z(\End(V)^H)=\bigoplus_jZ(\End(n_jV_j)^H)\simeq\R^{m-l}\oplus \C^l .
  $$
  
  \end{theorem}
  
   Let $\X^\infty_c(X/K)$ denote the Lie  algebra of compactly supported smooth strata preserving vector fields on $X/K$ (see \S\ref{quotientcomputations} for definitions).

  \begin{theorem}\label{quotienttheorem}
  Let $X$ and the $x_i$, $H_i$ and $W_i$ be as above. Then
  $$
  \H(\X^\infty_c(X/K))\simeq \bigoplus_i (Z(\End(W_i)^{H_i})/\lies_i)
  $$
  where each $\lie s_i$ is the Lie algebra of a torus $S_i$ lying in $Z(\End(W_i)^{H_i})$.
  \end{theorem}
  \noindent We will say more about the $S_i$ in \S \ref{quotientcomputations}.
  
\subsection{} This work  was done while attending the conference ``Diffeomorphisms and Related Fields'' held at Shinshu University, December 2005. The author thanks professors K.~Abe and K.~Fukui for the invitation and for their wonderful hospitality.
 
 \section{Vanishing of abelianizations}
\subsection{} In the following, let $\B^\infty_c(X)^K$ denote $[\X^\infty_c(X)^K,\X^\infty_c(X)^K]$ and let $\cc^\infty_c(X)^K$ denote the compactly supported smooth functions on $X$. Our first goal is to show that $\H(\X^\infty_c(X\times\R)^K)$ is zero.
 
 \begin{lemma} Let $A\in\xck$ and $B\in\xk$. Then $[A,B]\in\B^\infty_c(X)^K$.
  \end{lemma}
  
  \begin{proof}
  Let $g\in \cc^\infty_c(X)^K $ be  identically  $1$ on a neighborhood of $\supp A$.  Then $[A,gB]=g[A,B]+A(g)B=[A,B]\in\B^\infty_c(X)^K$.
  \end{proof}
  \begin{proposition}\label{crossrprop}
  Let $K$ act on $X\times\R$ with the given action on $X$ and the trivial action on $\R$. Then $\H(\X^\infty_c(X\times\R)^K)=0$.
  \end{proposition}
  \begin{proof} Let $t$ denote the usual coordinate function on $\R$ and let $g\in\cc^\infty_c(X\times\R)^K$. We show that $g\frac d{dt}\in\B^\infty_c(X\times\R)^K$.  For $x\in X$ and $s\in\R$ set $h(x,s)=\int_0^s g(x,u)\,du$. Then $h$ is smooth and $K$-invariant. Let $f\in\cc^\infty_c(X\times\R)^K$. Then
  $$
  [f\frac d{dt},h\frac d{dt}]=f\frac{dh}{dt}\frac d{dt}-h\frac{df}{dt}\frac d{dt} {\text{ and }}
  $$
  $$
  [\frac d{dt},fh\frac d{dt}]=f\frac {dh}{dt}\frac d{dt}+h\frac {df}{dt}\frac d{dt}.
  $$
  Hence $2fg\frac d{dt}\in\B^\infty_c(X\times\R)^K$. If $f$ equals $1/2$ on a neighborhood of $\supp g$, we obtain that $g\frac d{dt}\in\B^\infty_c(X\times\R)^K$.
  
  Now suppose that $A\in\X^\infty_c(X\times\R)^K$. By our result  above, we can assume that $A$ annihilates $t$. Set $B(x,s)=\int_0^s A(x,u)\,du$ and let $g\in\cc^\infty_c(X\times\R)^K$ equal $1$ on a neighborhood of $\supp A$. Then $[g\frac d{dt},B]=gA-B(g)\frac d{dt}$. We already know that $B(g)\frac d{dt}\in\B^\infty_c(X\times\R)^K$, hence $A\in\B^\infty_c(X\times\R)^K$. Thus $\H(\X^\infty_c(X\times\R)^K)=0$.
  \end{proof}
  
  \subsection{} \label{twistedproduct} Let $H$ be a closed subgroup of $K$ and $W$ an $H$-module. Then we have the twisted product $\twist KHW$ which is the quotient $(K\times W)/H$ where $h(k,w)=(kh\inv,hw)$, $h\in H$, $k\in K$ and $w\in W$. We denote the image of $(k,w)\in K\times W$ in $\twist KHW$ by $[k,w]$. Note that $\twist KHW$ is naturally a $K$-vector bundle and a real algebraic $K$-variety \cite{Algquotient}.
  
 Let $H\to\GL(W)$ be the slice representation at a point $x\in X$. By the differentiable slice theorem, a $K$-neighborhood of $Kx$ in $X$ is $K$-diffeomorphic to $\twist KHW$.  
  By Proposition \ref{crossrprop}, $\H(\X^\infty_c(\twist KHW)^K)=0$ if $W^H\neq (0)$. 
    
  Let $F$ be a closed $K$-stable subset of $X$. We say that {\it $\H(\X^\infty_c(X)^K)$ is supported on $F$\/} if $\H(\X^\infty_c(X\setminus F)^K)=0$. Using a partition of unity argument we can show
  \begin{corollary} \label{supportcor} Let $F=\{x\in X\mid W_x^{K_x}=0\}$. Then $\H(\xck)$ is supported on $F$.
  \end{corollary}
  
  \section{Local computations}
 \subsection{}\label{sumlocal} Our results above show that there is a discrete set of orbits $\{Kx_i\}$ such that 
  $$
  \H(\xck)\simeq\bigoplus_i\H(\X^\infty_c(K\ast^{H_i}W_i)^K)
  $$
  where $H_i=K_{x_i}$ and $W_i$ is the slice representation of $H_i$ at $x_i$. Thus it suffices to compute $\H(\X^\infty_c(K\ast^HV)^K)$ where $H$ is a closed subgroup of $K$, $V$ is an $H$-module   and $V^H=(0)$. This computation is the content  of the following theorem. 
 
  \begin{theorem}\label{localcomp} Let $H$ and $V$ be as above. Then
  $$
  \H(\X^\infty_c(K\ast^HV))\simeq\H(\liek^H/\lieh^H)\oplus\H(\End(V)^H).
  $$
  \end{theorem}

  \subsection{}\label{equations} Our proof of the theorem requires several lemmas. Set $Y:=K\ast^HV$. 
 Then
$$
  \X(Y)^K\simeq\X(K\times V)^{K\times H}/(\O(K\times V)\lieh)^{K\times H}
$$
(see  \cite[\S 4]{Lifting2})  where $H$ has the diagonal action (see \ref{twistedproduct}) on $K\times V$ (inducing an action of $\lieh$) and $\O(K\times V)$ denotes the polynomial functions on $K\times V$. Now
  $$
 \X(K\times V)^{K\times H}\simeq (\X(K)\otimes\O(V)\oplus\O(K)\otimes\X(V))^{K\times H} \simeq(\liek\otimes\O(V))^H\oplus(1\otimes\X(V)^H)  
   $$
  while
  $$
  (\O(K\times V)\lieh)^{K\times H}\simeq (\lieh\otimes\O(V))^H.
  $$
 
 \subsection{} We have the Euler operator $E\in\X(V)^H$, where if $x_1$, $x_2,\dots$ are coordinate functions on $V$, then $E=\sum_i x_i\frac{\pt}{\pt x_i}$. By the isomorphisms above,  $E$ can be considered as a $(K\times H)$-invariant vector field on $K\times V$ and as a $K$-invariant vector field on $Y$. 
 
  \begin{lemma}\label{elemma}
  Let $f\in\cc^\infty(Y)^K$. Then $f=E(h)$ for some $h\in\cc^\infty(Y)^K$ if and only if $f([e,0])=0$.
  \end{lemma}
  \begin{proof} Clearly the condition on $f$ is necessary. Suppose that $f([e,0])=0$. Since $f$ is $K$-invariant, it is determined by its restriction $g$ to $\{[e,v]\mid v\in V\}\simeq V$, where $g$ is $H$-invariant. Set $h(v)=\int_0^1 (1/t )g(tv)\,dt$. Then $h\in\cc^\infty(V)^H$ since $g(0)=0$. We have 
  $$
  E(h)(v)=\int_0^1  \frac 1t \sum_i x_i\frac {\pt g}{\pt x_i}(tv)t\,dt=\int_0^1\sum_i x_i \frac {\pt g}{\pt x_i}(tv)\,dt=\int_0^1\frac d{dt}g(tv)\,dt=g(v)-g(0)=g(v).
  $$
  \end{proof}
  \begin{corollary} \label{gEprop}Let $g\in\cc^\infty_c(Y)^K$ such that $g([e,0])=0$. Then $gE\in\B^\infty_c(Y)^K$.
  \end{corollary}
  \begin{proof} By Lemma \ref{elemma}, $g=E(h)$ for some $h\in\cc^\infty(Y)^K$. Let $f\in\cc^\infty_c(Y)^K$ such that $f$ is $1/2$ in a neighborhood of $\supp g$. Then, as in Proposition \ref{crossrprop}, 
  $$
  [E,fhE]+[fE,hE]=2fE(h)E=2fgE,
  $$
  so that $gE\in\B^\infty_c(Y)^K$.
  \end{proof}

 \subsection{} Since $Y$ is real algebraic, the results in \cite[\S 6]{Lifting} show that $\X^\infty(Y)\simeq\cc^\infty(Y)\otimes_{\O(Y)}\X(Y)$. For compactly supported sections we clearly have that $\X^\infty_c(Y)=\cc^\infty_c(Y)\X(Y)$.
  
\subsection{}   We  have an $E$-eigenspace decomposition
$$
\X(K\times V)^{K\times H}\simeq \bigoplus_{m\geq 0}(\liek\otimes\O(V)_m)^H\oplus(1\otimes\X(V)_m^H)
$$
 and similarly for $(\lieh\otimes\O(V))^H$. The weights that occur in $\X(V)^H$ are all positive since $V^H=(0)$. We have an induced  decomposition
  $$
  \X(Y)^K=\bigoplus_{m\geq 0}\X(Y)^K_m.
  $$
 \begin{remark}
 Since the sum only contains terms for $m\geq 0$, an element of $\X(Y)^K$ applied to an element of $\cc^\infty(Y)^K\simeq\cc^\infty(V)^H$ always vanishes at $[e,0]$.
 \end{remark}
  \begin{lemma} \label{fAomegaclemma} Let $A\in\X(Y)^K_m$ and let $f\in\cc^\infty_c(Y)^K$. Then $fA\in \B^\infty_c(Y)^K$ if 
  \begin{enumerate}
  \item  $m>0$ or
  \item $f([e,0])=0$.
  \end{enumerate}
  \end{lemma}
  \begin{proof} Suppose that $m>0$. Then $[(1/m)fE,A]=fA-(1/m)A(f)E$ where $A(f)E\in\B_c^\infty(Y)^K$ by Corollary \ref{gEprop}. Hence $fA\in\B^\infty_c(Y)^K$. If $m=0$ and $f([e,0])=0$, then let $h\in\cc^\infty(Y)^K$ be such that $E(h)=f$, and let $g\in\cc^\infty_c(Y)^K$. Then
  $$
  [gE,hA]=gE(h)A-hA(g)E=gfA-hA(g)E,
  $$
  where $hA(g)E\in\B^\infty_c(Y)^K$ by Corollary  \ref{gEprop}. We may arrange that $gfA=fA$, so $fA\in\B^\infty_c(Y)^K$.
  \end{proof}
  
  \begin{proof}[Proof of Theorem \ref{localcomp}]  We first define a map of Lie algebras $\phi\colon\X^\infty_c(Y)^K\to\X(Y)^K_0$. Let $B=\sum_{i=1}^m f_iB_i\in\X^\infty_c(Y)^K$ where $f_i\in\cc^\infty_c(Y)^K$ and $B_i\in\X(Y)^K_{m_i}$,  $i=1,\dots,m$. Define $\phi(B):=\sum_{m_i=0}f_i([e,0])B_i\in\X(Y)^K_0$.  It is obvious that $\phi$ is surjective. Suppose that $C$, $D\in\X(Y)^K$ are eigenvectors for $E$ and that $f$, $g\in\cc^\infty_c(Y)^K$.  
Then $[fC,gD]=fC(g)D-gD(f)C+fg[C,D]$ where $C(g)$ and $D(f)$ vanish at $[e,0]$. Thus $\phi([fC,gD])=(fg)(0)\phi([C,D])=(fg)(0)[\phi(C),\phi(D)]=[\phi(fC),\phi(gD)]$ . Now   $\phi$ induces $\tilde\phi\colon\H(\X^\infty_c(Y)^K)\to\H(\X(Y)^K_0)$, which is again surjective. Suppose that $B=\sum_i f_i B_i\in\Ker(\tilde\phi)$ where the $B_i$ are in $\X(Y)^K_0$. Then $\phi(B)=\sum_j[C_j,D_j]$ where   $C_j$, $D_j\in\X(Y)^K_0$ for all $j$. Let $f\in\cc^\infty_c(Y)^K$ such that $f$ is 1 on a neighborhood of $[e,0]$. Then $B -\sum_j[fC_j,fD_j]\in\B^\infty_c(Y)^K$. Hence $\tilde\phi$ is  an isomorphism. 
From our equations in \ref{equations} it follows that $\H(\X(Y)^K_0)\simeq\H(\liek^H/\lieh^H)\oplus\H(\End(V)^H)$.
  \end{proof}
  
  \begin{proof}[Proof of Theorem \ref{maintheorem}] The theorem is immediate from \ref{sumlocal} and Theorem \ref{localcomp}
  \end{proof}
  \begin{proof}[Proof of Theorem \ref{localtheorem}]
 Let $V=\oplus_{j=1}^m n_j V_j$ and $H$ be as in \ref{localtheorem}.    Then $\End(V)^H\simeq \oplus_j\End(n_jV_j)^H$. There are three cases to consider.
  
 Case 1: $\End(V_j)^H\simeq  \R$. Then $\End(n_jV_j)^H\simeq \liegl(n_j,\R)$ and $\H(\liegl(n_j,\R))\simeq Z(\liegl(n_j,\R))\simeq\R$.
  
  Case 2: $\End(V_j)^H\simeq\C$. Then $\End(n_jV_j)^H\simeq\liegl(n_j,\C)$ and $\H(\liegl(n_j,\C))\simeq Z(\liegl(n_j,\C))\simeq\C$.
  
 Case 3: $\End(V_j)^H\simeq\HH$, the quaternions. Then $\End(n_jV_j)^H\simeq\liegl(n_j,\HH)$ and we have that $\H(\liegl(n_j,\HH))\simeq Z(\liegl(n_j,\HH))\simeq\R$. The theorem follows.
  \end{proof}
  
    \section{Computations on the quotient}\label{quotientcomputations}
  We now consider the abelianization of the strata preserving vector fields on the quotient $X/K$. We recall a few facts about $X/K$ from \cite{Lifting}. Let $\pi\colon X\to X/K$ denote the canonical map, where $X/K$ is given the quotient topology. Then $X/K$ has a differentiable structure where for $U$ an open subset of $X/K$, $\cc^\infty(U)=\cc^\infty(\pi\inv(U))^K$. Let $H$ be a  closed subgroup   of $K$. Then we have the corresponding stratum $X^{(H)}:=\{x\in X\mid K_x$ is conjugate to $H\}$ and its image $(X/K)^{(H)}\subset X/K$. The {\it isotropy strata\/} $(X/K)^{(H)}\subset X/K$ and $X^{(H)}\subset X$ are  smooth and locally closed submanifolds and $\pi\colon X^{(H)}\to (X/K)^{(H)}$ is naturally a smooth fiber bundle (with structure group $N_K(H)/H$).  The number of isotropy  strata is locally finite  on $X$ and $X/K$. Let $\Der(\cc^\infty(X/K))$ denote the derivations of $\cc^\infty(X/K)$ and let $\X^\infty(X/K)$ denote those derivations that preserve the ideals of functions $I_{H_i}$ vanishing on   the isotropy  strata $(X/K)^{(H_i)}$ of $X/K$. Each element of $\X^\infty(X)^K$ restricts to a derivation of   $\cc^\infty(X/K)$, so there is a canonical map $\pi_*\colon \X^\infty(X)^K\to\Der(\cc^\infty(X/K))$. The main theorem of \cite{Lifting} is that $\Im\pi_*\subset\X^\infty(X/K)$ and that $\pi_*$ is surjective. Clearly $\pi_*$ is a homomorphism of Lie algebras so we have an induced surjection $\H(\X^\infty(X)^K)\to\H(\X^\infty(X/K))$. We only need to compute what happens in the case of $X=K\ast^HV$ where $H$ is a closed subgroup of $K$ and $V$ is an $H$-module such that $V^H=(0)$. 
  Let $V=\oplus_{j=1}^m n_j V_j$ as in Theorem \ref{localtheorem}.
  The following   has Theorem \ref{quotienttheorem} as a corollary.
  
  \begin{theorem}
  Assume that $\End(V_j)^H\simeq\C$ if and only if $j\leq l$ where $l\leq m$. Let $T$ be the corresponding torus $(S^1)^l\subset \prod_{j=1}^l Z(\End(V_j)^H)$. Then $T$ acts on $V$ commuting with the action of $H$, and we have an induced map $T\to\Aut(V/H)$. Let $S$ denote the kernel    where $\dim S=k$.   Then
    $$
  \H(\X_c^\infty((K\ast^HV)/K))\simeq\H(\X(V/H))\simeq \R^{m-l+k}\oplus\C^{l-k}.
  $$
  \end{theorem}
  \begin{proof}
  We have the canonical surjection of Lie algebras $\pi_*\colon\End(V)^H\to\X_0(V/H)$ and  $\pi_*$ induces a surjection of $\H(\End(V)^H)$ onto $\H(\X(V/H))$.
For every $j$ we have the identity $\Id_j\in\End(n_jV_j)^H$  and clearly these elements give linearly independent derivations of $\O(V)^H$. Now consider the action of $T$ on $V/H$ and its kernel $S$. Then $\lie s$ is the kernel of the restriction of $\pi_*$  to the center of $\End(V)^H$, so that $\lie s$ is the kernel on homology.  
   \end{proof}
   \begin{example}
   Suppose that $H$ is a torus acting faithfully on $V$ and  $V=\sum_{j=1}^m n_jV_j$ where $V^H=(0)$  as in Theorem \ref{localtheorem}. Then $\lies\simeq \lieh$ and $\H(\X(V/H))\simeq \R^k\oplus\C^{m-k}$ where $k=\dim H$.  
   \end{example}
  \begin{example}
  Let $V=\C^n\oplus\wedge^2\C^n$ with the canonical action of $\SU(n,\C)$, $n\geq 3$. Then $T$ has dimension 2 and $S$ has dimension 1. See \cite[Table I]{Lifting}.
    \end{example}


\newcommand{\noopsort}[1]{} \newcommand{\printfirst}[2]{#1}
\newcommand{\singleletter}[1]{#1} \newcommand{\switchargs}[2]{#2#1}
\providecommand{\bysame}{\leavevmode\hbox to3em{\hrulefill}\thinspace}

\end{document}